\documentclass[a4paper,10pt,leqno,twoside]{article}

\usepackage{amsmath, amsthm, amssymb,bm}

\numberwithin{equation}{section}

\newcommand{\bC}{\mathbb{C}}
\newcommand{\bR}{\mathbb{R}}
\newcommand{\bQ}{\mathbb{Q}}
\newcommand{\bT}{\mathbb{T}}

\newcommand{\bZ}{\mathbb{Z}}

\newcommand{\mf}[1]{\mathfrak{#1}}
\newcommand{\mr}[1]{\mathrm{#1}}
\newcommand{\mcal}[1]{\mathcal{#1}}

\def\trans#1{#1^{\mathrm{T}}} 
\def\dd{\mathrm{d}}

\def\bx{\bm{x}}
\def\by{\bm{y}}
\def\bz{\bm{z}}

\newcommand{\CBE}[1]{\mathrm{C}\beta\mathrm{E}_{#1}}

\def\({ \left( }
\def\){ \right)}
\def\[{ \left[ }
\def\]{ \right]}
\theoremstyle{plain}
\newtheorem{thm}{Theorem}[section]
\newtheorem{prop}[thm]{Proposition}

\newtheorem{cor}[thm]{Corollary}
\theoremstyle{definition}

\newtheorem{example}{Example}[section]
\newtheorem{remark}{Remark}[section]
\theoremstyle{conjecture}

\theoremstyle{problem}

\title{\bfseries  Moments of characteristic polynomials for
compact symmetric spaces
and Jack polynomials}
\author{\textsc{Sho MATSUMOTO}
\thanks{Research Fellow of the Japan Society 
for the Promotion of Science, partially supported by Grant-in-Aid 
for Scientific Research (C) No. 17006193.}%
}
\date{}

\begin{document}
%>>>>>>>>>>>>>>>>>>>>>>>>>>>>>>>>>>>>>>>>>
%\setlength{\baselineskip}{14pt}
\maketitle

\begin{abstract}
We express the averages of products of characteristic polynomials for random matrix ensembles associated with compact symmetric spaces in terms of
Jack polynomials or 
Heckman and Opdam's Jacobi polynomials depending on the root system of the space.
We also give explicit expressions for 
the asymptotic behavior of these averages in the limit as the matrix size goes to infinity.

\noindent
{\bf MSC-class}: primary 15A52; secondary 33C52, 05E05. \\
\noindent
{\bf Keywords}: characteristic polynomial, random matrix, Jacobi polynomial, Jack polynomial,
Macdonald polynomial, compact symmetric space.
\end{abstract}

%((((((((((((((    Main Text follows    ))))))))))))))))

%
%%%%%%%%%%%%%%%%%%%%%%%%%%%%%%%%%%%%%%%%%%%%%%%%%%%%%%%%%%%%%%%%%%%%%%%%%%%%%%%%%%%%%%%%%%%%%%%%
\section{Introduction}
%%%%%%%%%%%%%%%%%%%%%%%%%%%%%%%%%%%%%%%%%%%%%%%%%%%%%%%%%%%%%%%%%%%%%%%%%%%%%%%%%%%%%%%%%%%%%%%%
%

In recent years, there has been considerable interest in the averages 
of the characteristic polynomials of random matrices. 
This work is motivated by the connection with Riemann zeta functions and $L$-functions 
identified by Keating and Snaith \cite{KS_zetafunctions, KS_Lfunctions}.
The averages of the characteristic polynomials 
in the cases of compact classical groups and Hermitian matrix ensembles
have already calculated,
see \cite{Mehta} and references in \cite{BG}.
In these studies, Bump and Gamburd \cite{BG} obtain simple proofs 
for the cases corresponding to compact classical groups
by using symmetric polynomial theory.
Our aim in this note is to use their technique to calculate averages of the characteristic polynomials
for random matrix ensembles 
associated with compact symmetric spaces.

We deal with the compact symmetric spaces $G/K$ classified by Cartan,
where $G$ is a compact subgroup in $GL(N,\bC)$ for some positive integer $N$,
and $K$ is a closed subgroup of $G$.
Assume $G/K$ is realized as a subspace $S$ in $G$, i.e., $S \simeq G/K$,
and the probability measure $\dd M$ on $S$ is then induced from $G/K$.
We call the probability space $(S, \dd M)$ the random matrix ensemble
associated with $G/K$.

For example, $U(n)/O(n)$ is the symmetric space with a restricted root system of type A,
and is realized by $S=\{M \in U(n) \ | \ M = \trans{M} \}$.
Here $\trans{M}$ stands for the transposed matrix of $M$ 
while $U(n)$ and $O(n)$ denote the unitary and orthogonal group 
of matrices or order $n$ respectively.
The induced measure $\dd M$ on $S$ satisfies the invariance 
$\dd (H M \trans{H})= \dd M$ for any $H \in U(n)$.
This random matrix ensemble $(S, \dd M)$ is well known 
as the circular orthogonal ensemble (COE for short), see e.g. \cite{Dyson,Mehta}.

We also  consider the classical compact Lie groups $U(n)$, $SO(n)$, and $Sp(2n)$.
Regarding these groups as symmetric spaces,
 the random matrix space $S$ is just the group itself with its Haar measure.
 
The compact symmetric spaces 
studied by Cartan are divided into A and BC type
main branches according to their root systems.
There are three symmetric spaces of type A, with  
their corresponding matrix ensembles called 
circular orthogonal, unitary, and symplectic ensembles.  
For these ensembles,
the probability density functions (p.d.f.) for the eigenvalues 
are proportional to
$$
\Delta^{\mr{Jack}}(\bz;2/\beta)= \prod_{1 \le i<j \le n} |z_i -z_j|^{\beta}, 
$$
with $\beta=1,2,4$,
where $\bz =(z_1,\dots, z_n)$, with $|z_i|=1$, denotes the sequence of eigenvalues 
of the random matrix.
We will express the average of the product of
characteristic polynomials $\det(I+ xM)$ for a random matrix $M$
 as a Jack polynomial
(\cite[Chapter VI-10]{Mac})
of a rectangular-shaped Young diagram.
Jack polynomials are orthogonal with respect to the weight function 
$\Delta^{\mr{Jack}}$.
Our theorems are obtained in a simple algebraic way, and 
contain results given in \cite{KS_zetafunctions}.

For compact symmetric spaces of type BC root systems, 
the corresponding p.d.f. is given by
$$
\Delta^{\mr{HO}}(\bz;k_1,k_2,k_3) =
\prod_{1 \le i <j \le n} |1-z_i z_j^{-1}|^{2k_3} |1-z_i z_j|^{2k_3}
\cdot \prod_{1 \le j \le n} |1-z_j|^{2k_1} |1-z_j^2|^{2k_2}.
$$
Here 
the $k_i$'s denote multiplicities of roots in 
the root systems of the symmetric spaces.
For example, the p.d.f. induced from the symmetric space $SO(4n+2)/(SO(4n+2) \cap Sp(4n+2))$
is proportional to $\Delta^{\mr{HO}}(\bz;2, \frac{1}{2},2)$. 
For this class of compact symmetric spaces,
Opdam and Heckman's Jacobi polynomials (\cite{Diejen, Heckman}),
which are orthogonal with respect to $\Delta^{\mr{HO}}$,
will play the same role as Jack polynomials for type A cases.
Namely, we will express the average of the product of
characteristic polynomials $\det(I+ xM)$ as the Jacobi polynomial
of a rectangular-shaped diagram.

This paper is organized as follows:

Our main results, which are expressions for the averages of products of characteristic polynomials,
will be given \S 6.
As described above, the symmetric spaces corresponding to the two root systems, type A and BC, will be discussed separately.
For  type A spaces, we use Jack polynomial theory. 
These discussions can be generalized to Macdonald polynomials.
Thus, after preparations in \S 2,
we give some generalized identities 
involving Macdonald polynomials and a generalization of the weight function $\Delta^{\mr{Jack}}$
in \S 3 and \S 4.
In particular, we obtain $q$-analogues of Keating and Snaith's formulas \cite{KS_zetafunctions}
for the moments of characteristic polynomials
and a generalization of the strong Szeg\"{o} limit theorem for Toeplitz determinants.
These identities are reduced to 
characteristic polynomial expressions 
for symmetric spaces of the A type root system in \S 6.1 - \S 6.3.
On the other hand, for type BC spaces, we employ Opdam and Heckman's Jacobi polynomials.
We review the definition and several properties of these polynomials in \S 5, 
while in \S 6.4 - \S 6.12 we apply  them to obtain expressions 
for the products of characteristic polynomials of random matrix ensembles 
associated with symmetric spaces of type BC.

%
%%%%%%%%%%%%%%%%%%%%%%%%%%%%%%%%%%%%%%%%%%%%%%%%%%%%%%%%%%%%%%%%%%%%%%%%%%%%%%%%%%%%%%%%%%%%%%%%
\section{Basic Properties of  Macdonald symmetric functions}
%%%%%%%%%%%%%%%%%%%%%%%%%%%%%%%%%%%%%%%%%%%%%%%%%%%%%%%%%%%%%%%%%%%%%%%%%%%%%%%%%%%%%%%%%%%%%%%%
%

We recall the definition of Macdonald symmetric functions, see \cite[Chapter VI]{Mac} for details.
Let $\lambda$ be a partition, i.e., $\lambda=(\lambda_1,\lambda_2,\dots)$
is a weakly decreasing ordered sequence of non-negative integers with finitely many non-zero entries.
Denote by $\ell(\lambda)$ the number of non-zero $\lambda_j$ and 
by $|\lambda|$ the sum of all $\lambda_j$.
These values $\ell(\lambda)$ and $|\lambda|$ are called
the length and weight of $\lambda$ respectively.
We identify $\lambda$ with the associated Young diagram 
$\{(i,j) \in \bZ^2 \ | \ 1 \le j \le \lambda_i \}$.
The conjugate partition $\lambda'=(\lambda'_1,\lambda'_2,\dots)$ 
is determined by the transpose of the Young diagram $\lambda$.
It is sometimes convenient to write this partition in the form $\lambda=(1^{m_1} 2^{m_2} \cdots )$,
where $m_i=m_i(\lambda)$ is the multiplicity of $i$ in $\lambda$
and is given by $m_i=\lambda'_i-\lambda'_{i+1}$.
For two partitions $\lambda$ and $\mu$,
we write $\lambda \subset \mu$  
if $\lambda_i \le \mu_i$ for all $i$.
In particular, the notation $\lambda \subset (m^n)$
means that 
$\lambda$ satisfies $\lambda_1 \le m$ and $\lambda_1' \le n$.
The dominance ordering associated with the root system of type A is defined as follows:
for two partitions $\lambda=(\lambda_1,\lambda_2,\dots)$ and $\mu=(\mu_1,\mu_2,\dots)$,
$$
\mu \le_{\mr{A}} \lambda \qquad \Leftrightarrow \qquad 
|\lambda|=|\mu|
\quad \text{and} \quad
\mu_1 + \cdots+\mu_i \le \lambda_1+ \cdots +\lambda_i \quad \text{for all $i \ge 1$}.
$$

Let $q$ and $t$ be real numbers such that both $|q|<1$ and $|t|<1$. 
Put $F=\bQ(q,t)$ and $\bT^n=\{\bz =(z_1,\dots,z_n) \ | \ |z_i|=1 \ (1 \le i \le n)\}$.
Denote by $F[x_1,\dots,x_n]^{\mf{S}_n}$ the algebra of 
symmetric polynomials in variables $x_1,\dots,x_n$.
Define an inner product on $F[x_1,\dots,x_n]^{\mf{S}_n}$
by
$$
\langle f, g \rangle_{\Delta^{\mr{Mac}}} = \frac{1}{n!}
\int_{\bT^n} f(\bz) g(\bz^{-1}) \Delta^{\mr{Mac}}(\bz;q,t) \dd \bz
$$
with
$$
\Delta^{\mr{Mac}}(\bz;q,t)= \prod_{1 \le i<j \le n} \Bigg|
\frac{(z_i z_j^{-1};q)_\infty}{(t z_i z_j^{-1};q)_\infty} \Bigg|^2,
$$
where $\bz^{-1}=(z_1^{-1},\dots,z_n^{-1})$ and $(a;q)_\infty= \prod_{r=0}^\infty(1-aq^r)$.
Here $\dd \bz$ is the normalized Haar measure on $\bT^n$.

For a partition $\lambda$ of length $\ell(\lambda) \le n$, put 
\begin{equation} \label{eq:monomialA}
m_{\lambda}^{\mr{A}} (x_1,\dots,x_n) = 
\sum_{\nu=(\nu_1,\dots,\nu_n) \in \mf{S}_n \lambda} x_1^{\nu_1} \cdots x_n^{\nu_n},
\end{equation}
where the sum runs over the $\mf{S}_n$-orbit 
$\mf{S}_n \lambda = \{ (\lambda_{\sigma(1)},\dots, \lambda_{\sigma(n)}) \ | \ \sigma \in \mf{S}_n\}$.
Here we add the suffix ``A'' because $\mf{S}_n$ is the Weyl group of type A. 
Then 
Macdonald polynomials (of type A)
$P_\lambda^{\mr{Mac}}=P_{\lambda}^{\mr{Mac}}(x_1,\dots,x_n;q,t)
\in F[x_1,\dots,x_n]^{\mf{S}_n}$ are characterized
by the following conditions:
$$
P_{\lambda}^{\mr{Mac}} = m_{\lambda}^{\mr{A}} + \sum_{\mu <_{\mr{A}} \lambda} u_{\lambda \mu} 
m_{\mu}^{\mr{A}} 
\quad \text{with $u_{\lambda\mu} \in F$}, \qquad\qquad
\langle P_{\lambda}^{\mr{Mac}}, P_{\mu}^{\mr{Mac}} 
\rangle_{\Delta^{\mr{Mac}}}=0 \quad \text{if $\lambda \not=\mu$}.
$$

Denote by $\Lambda_F$ the $F$-algebra of symmetric functions in infinitely many variables 
$\bx=(x_1,x_2,\dots)$.
That is, an element $f =f(\bx) \in \Lambda_F$ is determined by
the sequence $(f_n)_{n \ge 0}$ of
polynomials $f_n$ in $F[x_1,\dots,x_n]^{\mf{S}_n}$, where these polynomials satisfy $\sup_{n \ge 0} \deg (f_n) < \infty$ and 
$f_m(x_1,\dots,x_n,0,\dots,0)=f_n(x_1,\dots,x_n)$ for any $m \ge n$,
see \cite[Chapter I-2]{Mac}.
Macdonald polynomials satisfy the stability property 
$$
P^{\mr{Mac}}_\lambda(x_1,\dots,x_n,x_{n+1};q,t) \Big|_{x_{n+1}=0} =
P^{\mr{Mac}}_\lambda(x_1,\dots,x_n;q,t)
$$
for any partition $\lambda$ of length
$\ell(\lambda) \le n$,
and therefore for all partitions $\lambda$, 
{\it Macdonald functions} $P_{\lambda}^{\mr{Mac}}(\bx ;q,t)$
can be defined.

For each square $s=(i,j)$ of the diagram $\lambda$, let
$$
a(s)=\lambda_i-j, \qquad a'(s)=j-1, \qquad l(s)= \lambda'_j-i, \qquad l'(s)= i-1.
$$
These  numbers are called the arm-length, arm-colength, leg-length, and leg-colength respectively.
Put
$$
c_{\lambda}(q,t)= \prod_{s \in \lambda} (1-q^{a(s)}t^{l(s)+1}), \qquad
c_{\lambda}'(q,t)= \prod_{s \in \lambda} (1-q^{a(s)+1} t^{l(s)}).
$$
Note that $c_{\lambda}(q,t)= c'_{\lambda'}(t,q)$.
Defining the $Q$-function by 
$Q_{\lambda}(\bx;q,t)= c_{\lambda}(q,t) c'_{\lambda}(q,t)^{-1} P_\lambda(\bx;q,t)$,
we have the dual Cauchy identity \cite[Chapter VI (5.4)]{Mac}
\begin{align}
& \sum_{\lambda} P_\lambda(\bx;q,t) P_{\lambda'}(\by;t,q)=
\sum_{\lambda} Q_\lambda(\bx;q,t) Q_{\lambda'}(\by;t,q) \label{EqDualCauchy} \\
=& \prod_{i\ge 1} \prod_{j\ge 1} (1+x_i y_j)
=\exp \(\sum_{k=1}^\infty \frac{(-1)^{k-1}}{k}p_k(\bx)p_k(\by) \), \notag 
\end{align}
where $\by=(y_1,y_2,\dots)$.
Here $p_k$ is the power-sum function $p_k(\bx)=x_1^k+x_2^k+ \cdots$.

We define the generalized factorial $(a)_{\lambda}^{(q,t)}$ by
$$
(a)_{\lambda}^{(q,t)} = \prod_{s \in \lambda} (t^{l'(s)} - q^{a'(s)} a).
$$
Let $u$ be an indeterminate and define the homomorphism $\epsilon_{u,t}$ from $\Lambda_F$
to $F$ by
\begin{equation} \label{EqSpecialPowerSum}
\epsilon_{u,t}(p_r) = \frac{1-u^r}{1-t^r} \qquad \text{for all $r \ge 1$}.
\end{equation}
In particular, we have $\epsilon_{t^n,t}(f)= f(1,t,t^2,\dots, t^{n-1})$ 
for any $f \in \Lambda_F$.
Then we have (\cite[Chapter VI (6.17)]{Mac})
\begin{equation} \label{EqSpecialMac}
\epsilon_{u,t}(P_{\lambda}^{\mr{Mac}})= \frac{(u)_{\lambda}^{(q,t)}}{c_{\lambda}(q,t)}.
\end{equation}
Finally, the following orthogonality property is satisfied for any two partitions $\lambda$ and $\mu$ of length $\le n$:
\begin{equation} \label{EqOrthogonality}
\langle P_{\lambda}^{\mr{Mac}}, Q_{\mu}^{\mr{Mac}} \rangle_{\Delta^{\mr{Mac}}}= 
\delta_{\lambda \mu} \langle 1,1 \rangle_{\Delta^{\mr{Mac}}}
\prod_{s \in \lambda} 
\frac{1-q^{a'(s)}t^{n-l'(s)}}{1-q^{a'(s)+1}t^{n-l'(s)-1}}.
\end{equation}

%
%%%%%%%%%%%%%%%%%%%%%%%%%%%%%%%%%%%%%%%%%%%%%%%%%%%%%%%%%%%%%%%%%%%%%%%%%%%%%%%%%%%%%%%%%%%%%%%%
\section{Averages with respect to $\Delta^{\mr{Mac}}(\bz;q,t)$} \label{sectionMacAverage}
%%%%%%%%%%%%%%%%%%%%%%%%%%%%%%%%%%%%%%%%%%%%%%%%%%%%%%%%%%%%%%%%%%%%%%%%%%%%%%%%%%%%%%%%%%%%%%%%
%

As in the previous section, we assume
$q$ and $t$ are real numbers in the interval $(-1,1)$.
For a Laurent polynomial $f$ in variables $z_1,\dots,z_n$, we define
$$
\langle f \rangle_{n}^{(q,t)} = 
\frac{\int_{\bT^n} f(\bz) \Delta^{\mr{Mac}}(\bz;q,t) \dd \bz}
{\int_{\bT^n} \Delta^{\mr{Mac}}(\bz;q,t) \dd \bz}.
$$
In this section,
we calculate averages of the products of the polynomial
$$
\Psi^{\mr{A}}(\bz;\eta)= \prod_{j=1}^n (1+ \eta z_j), \qquad \eta \in \bC
$$
with respect to $\langle \cdot \rangle_{n}^{(q,t)}$.
Denoting the eigenvalues of a unitary matrix $M$ by $z_1,\dots,z_n$, the polynomial $\Psi^{\mr{A}}(\bz;\eta)$ is the characteristic polynomial $\det(I+\eta M)$. 
 
The following theorems will induce averages of the products of characteristic polynomials
for random matrix ensembles associated with   type A root systems,
see \S \ref{sectionCBEq} and \S \ref{subsectionA} - \S \ref{subsectionAII} below. 

\begin{thm} \label{ThmAverageMac}
Let  $K$ and $L$ be positive integers.
Let $\eta_1,\dots, \eta_{L+K}$ be
complex numbers such that $\eta_j \not=0 \ ( 1\le j \le L)$.
Then we have 
$$
\left\langle \prod_{l=1}^L \Psi^{\mr{A}}(\bz^{-1};\eta_l^{-1}) \cdot \prod_{k=1}^K
\Psi^{\mr{A}}(\bz;\eta_{L+k}) \right\rangle_{n}^{(q,t)} =
(\eta_1 \cdots \eta_L)^{-n}
\cdot P_{(n^L)}^{\mr{Mac}} (\eta_1, \dots, \eta_{L+K};t,q).
$$
\end{thm}

\begin{proof}
By the dual Cauchy identity \eqref{EqDualCauchy}, we have
\begin{align*}
& \prod_{l=1}^L
\Psi^{\mr{A}}(\bz^{-1};\eta_l^{-1}) \cdot \prod_{k=1}^K
\Psi^{\mr{A}}(\bz;\eta_{L+k})
= \prod_{l=1}^L \eta_l^{-n} \cdot (z_1 \cdots z_n)^{-L} \cdot 
\prod_{k=1}^{L+K} \prod_{j=1}^n (1+\eta_k z_j) \\
=& \prod_{l=1}^L \eta_l^{-n} \cdot (z_1 \cdots z_n)^{-L}
\sum_{\lambda} Q_{\lambda}^{\mr{Mac}}(\eta_1,\dots,\eta_{L+K};t,q) 
Q_{\lambda'}^{\mr{Mac}}(\bz;q,t).
\end{align*}
Therefore, since $P_{(L^n)}^{\mr{Mac}}(\bz;q,t)=(z_1\cdots z_n)^L$ (\cite[Chapter VI (4.17)]{Mac}),
we see that
\begin{align*}
 \left\langle \prod_{l=1}^L \Psi^{\mr{A}}(\bz^{-1};\eta_l^{-1}) \cdot \prod_{k=1}^K
\Psi^{\mr{A}}(\bz;\eta_{L+k}) \right\rangle_{n}^{(q,t)} 
=& \prod_{l=1}^L \eta_l^{-n} 
\sum_{\lambda} Q_{\lambda}^{\mr{Mac}}(\eta_1,\dots,\eta_{L+K};t,q) 
\frac{\langle Q_{\lambda'}^{\mr{Mac}},
P_{(L^n)}^{\mr{Mac}} \rangle_{ \Delta^{\mr{Mac}} }   }
{\langle 1, 1 \rangle_{\Delta^{\mr{Mac}} } }
 \\
=& \prod_{l=1}^L \eta_l^{-n} \cdot Q_{(n^L)}^{\mr{Mac}}(\eta_1,\dots,\eta_{L+K};t,q)
\prod_{s \in (L^n)} 
\frac{1-q^{a'(s)}t^{n-l'(s)}}{1-q^{a'(s)+1}t^{n-l'(s)-1}}
\end{align*}
by the orthogonality property \eqref{EqOrthogonality}.
It is easy to check that
$$
\prod_{s \in (L^n)} 
\frac{1-q^{a'(s)}t^{n-l'(s)}}{1-q^{a'(s)+1}t^{n-l'(s)-1}} 
=\frac{c_{(L^n)}(q,t)}{c'_{(L^n)}(q,t)}= \frac{c'_{(n^L)}(t,q)}{c_{(n^L)}(t,q)},
$$
and so we obtain the claim.
\end{proof}

It may be noted that the present proof of Theorem \ref{ThmAverageMac} is similar to the 
corresponding one in \cite{BG}.

\begin{cor} \label{CorMomentValue}
For each positive integer $k$ and $\xi \in \bT$,
we have
$$
\left\langle \prod_{i=0}^{k-1} | \Psi^{\mr{A}}(\bz; q^{i+1/2} \xi )|^2 \right\rangle_{n}^{(q,t)}
=
\prod_{i=0}^{k-1} \prod_{j=0}^{n-1} \frac{1-q^{k+i+1} t^j}{1-q^{i+1} t^j}.
$$
\end{cor}

\begin{proof}
Set $L=K=k$ and $\overline{\eta_i}^{-1} =\eta_{i+k} =q^{i-1/2} \xi \ (1 \le i \le k)$ 
in Theorem \ref{ThmAverageMac}.
Then we have
\begin{align*}
\left\langle \prod_{i=0}^{k-1} | 
\Psi^{\mr{A}}(\bz;q^{i+1/2} \xi)|^2 \right\rangle_{n}^{(q,t)}
=& \prod_{i=0}^{k-1} q^{(i+1/2)n} \cdot P_{(n^k)} (q^{-k+1/2}, q^{-k+3/2},  \dots, q^{-1/2},
q^{1/2}, \cdots, q^{k-1/2};t,q)  \\
=& q^{n k^2/2} \cdot q^{(-k+1/2)kn} P_{(n^k)}(1,q,\cdots, q^{2k-1};t,q)  \\
=& q^{-n k(k-1)/2} \epsilon_{q^{2k},q} (P_{(n^k)}(\cdot;t,q)). 
\end{align*}
From expression \eqref{EqSpecialMac}, the right-hand side of the above expression equals
$$
 q^{-n k(k-1)/2} \frac{(q^{2k})_{(n^k)}^{(t,q)}}{c_{(n^k)}(t,q)}
=q^{-n k(k-1)/2} \prod_{i=1}^k \prod_{j=1}^n \frac{q^{i-1}-t^{j-1} q^{2k}}{1-t^{n-j}q^{k-i+1}} 
=  \prod_{j=0}^{n-1} \prod_{i=1}^k \frac{1-t^{j} q^{2k-i+1}}{1-t^j q^{k-i+1}},
$$
and the result follows.
\end{proof}

Kaneko \cite{Kaneko2} defines
the multivariable $q$-hypergeometric function associated with Macdonald polynomials by
$$
{_2 \Phi_1}^{(q,t)}(a,b;c;x_1,\dots,x_n)= \sum_\lambda
\frac{ (a)_{\lambda}^{(q,t)} (b)_{\lambda}^{(q,t)}}{(c)_{\lambda}^{(q,t)}}
\frac{P^{\mr{Mac}}_{\lambda}(x_1,\dots,x_n;q,t)}{c'_{\lambda}(q,t)},
$$
where $\lambda$ runs over all partitions of length $\ell(\lambda) \le n$.
The $q$-shifted moment 
$\left\langle \prod_{i=0}^{k-1} | \Psi^{\mr{A}}(\bz;q^{i+1/2} \xi )|^2 \right\rangle_{n}^{(q,t)}$
given 
in Corollary \ref{CorMomentValue}
can also be expressed as a special value of the generalized
$q$-hypergeometric function ${_2 \Phi_1}^{(q,t)}$
as follows:

\begin{prop} \label{PropMomentHypergeometric}
For any complex number with $|\eta|<1$ and real number $u$,
$$
\left\langle \prod_{j=1}^n \left|
 \frac{(\eta z_j;q)_\infty}{(\eta z_j u; q)_{\infty}}  \right|^2 \right\rangle_{n}^{(q,t)}
= {_2 \Phi_1}^{(q,t)}(u^{-1},u^{-1} ;q t^{n-1};
(u|\eta|)^2, (u|\eta|)^2 t,\dots, (u|\eta|)^2t^{n-1}).
$$
In particular, letting $u=q^k$ and $\eta=q^{1/2}\xi$ with $\xi \in \bT$,
we have
$$
\left\langle \prod_{i=0}^{k-1} | \Psi^{\mr{A}}(\bz;q^{i+1/2} \xi)|^2 \right\rangle_{n}^{(q,t)}
= {_2 \Phi_1}^{(q,t)}(q^{-k},q^{-k} ;q t^{n-1};
q^{2k+1}, q^{2k+1}t , \dots, q^{2k+1} t^{n-1}).
$$
\end{prop}

\begin{proof}
A simple calculation gives
$$
\prod_{j=1}^n
 \frac{(\eta z_j;q)_\infty}{(\eta z_j u; q)_{\infty}} 
= \exp \(\sum_{k=1}^\infty \frac{(-1)^{k-1}}{k} \frac{1-u^k}{1-q^k} 
p_{k}(-\eta z_1, \dots, -\eta z_n) \).
$$
From expressions \eqref{EqDualCauchy} and \eqref{EqSpecialPowerSum}, we have
$$
\prod_{j=1}^n
 \frac{(\eta z_j;q)_\infty}{(\eta z_j u; q)_{\infty}} 
= \sum_{\lambda} (-\eta)^{|\lambda|} 
\epsilon_{u,q}(Q^{\mr{Mac}}_{\lambda'}(\cdot;t,q)) Q^{\mr{Mac}}_{\lambda}(\bz;q,t)
= \sum_{\lambda} (-\eta)^{|\lambda|} 
\epsilon_{u,q}(P^{\mr{Mac}}_{\lambda'}(\cdot;t,q)) P^{\mr{Mac}}_{\lambda}(\bz;q,t).
$$
Thus we have
$$
\prod_{j=1}^n \left|
 \frac{(\eta z_j;q)_\infty}{(\eta z_j u; q)_{\infty}}  \right|^2
= \sum_{\lambda, \mu} (-\eta)^{|\lambda|} (-\overline{\eta})^{|\mu|}  
\epsilon_{u,q}(P^{\mr{Mac}}_{\lambda'}(\cdot;t,q)) \epsilon_{u,q}(Q^{\mr{Mac}}_{\mu'}(\cdot;t,q))
 P^{\mr{Mac}}_{\lambda}(\bz;q,t) Q^{\mr{Mac}}_{\mu}(\bz^{-1};q,t).
$$
The average is given by
\begin{align*}
\left\langle \prod_{j=1}^n \left|
 \frac{(\eta z_j;q)_\infty}{(\eta z_j u; q)_{\infty}}  \right|^2 \right\rangle_{n}^{(q,t)} 
=& \sum_{\lambda} |\eta|^{2|\lambda|} 
\epsilon_{u,q}(P^{\mr{Mac}}_{\lambda'}(\cdot;t,q)) 
\epsilon_{u,q}(Q^{\mr{Mac}}_{\lambda'}(\cdot;t,q))
\frac{ \langle  P^{\mr{Mac}}_{\lambda}, Q^{\mr{Mac}}_{\lambda}
\rangle_{\Delta^{\mr{Mac}}} }
{\langle  1, 1 \rangle_{\Delta^{\mr{Mac}}} } \\
=& \sum_{\lambda} |\eta|^{2|\lambda|} 
\frac{\{(u)_{\lambda'}^{(t,q)}\}^2}{c_{\lambda'}(t,q) c'_{\lambda'}(t,q)}
\prod_{s \in \lambda} \frac{1-q^{a'(s)}t^{n-l'(s)}}{1-q^{a'(s)+1}t^{n-l'(s)-1}}
\end{align*}
by expression \eqref{EqSpecialMac} and the orthogonality property \eqref{EqOrthogonality}.
It is easy to check that
\begin{align*}
&(u)_{\lambda'}^{(t,q)} =(-u)^{|\lambda|} (u^{-1})_{\lambda}^{(q,t)}, \qquad
c_{\lambda'}(t,q) c'_{\lambda'}(t,q)= c_{\lambda}(q,t) c'_{\lambda}(q,t), \\
&\prod_{s \in \lambda} \frac{1-q^{a'(s)}t^{n-l'(s)}}{1-q^{a'(s)+1}t^{n-l'(s)-1}}
= \prod_{s \in \lambda} \frac{t^{l'(s)}-q^{a'(s)}t^{n}}{t^{l'(s)} -q^{a'(s)+1}t^{n-1}}
= \frac{(t^n)_{\lambda}^{(q,t)}}{(q t^{n-1})_{\lambda}^{(q,t)}}.
\end{align*}
Finally, we obtain
$$
\left\langle \prod_{j=1}^n \left|
 \frac{(\eta z_j;q)_\infty}{(\eta z_j u; q)_{\infty}}  \right|^2 \right\rangle_{n}^{(q,t)} 
= \sum_{\lambda} (u|\eta|)^{2|\lambda|} \frac{\{(u^{-1})_{\lambda}^{(q,t)}\}^2}{(q t^{n-1})_{\lambda}^{(q,t)}}
\frac{P^{\mr{Mac}}_{\lambda}(1,t,\dots,t^{n-1};q,t)}{c'_{\lambda}(q,t)},
$$
which equals ${_2 \Phi_1}^{(q,t)}(u^{-1},u^{-1} ;q t^{n-1};
(u|\eta|)^2,\dots, (u|\eta|)^2t^{n-1})$.
\end{proof}

Now we derive the asymptotic behavior of the moment of $|\Psi(\bz;\eta)|$
when $|\eta| < 1$ in the limit as $n \to \infty$.
The following theorem is a generalization of the well-known
strong Szeg\"{o} limit theorem as stated in \S \ref{subsectionCBEJack} below.

\begin{thm} \label{Thm:SzegoMacdonald}
Let $\phi(z)=\exp(\sum_{k \in \bZ} c(k) z^k)$ be a function on $\bT$ and assume
\begin{equation} \label{Eq:AssumptionSzego}
\sum_{k \in \bZ} |c(k)|< \infty \qquad \text{and} \qquad 
\sum_{k \in \bZ} |k| |c(k)|^2 < \infty.
\end{equation}
Then we have
$$
\lim_{n \to \infty} e^{-n c(0)} 
\left\langle \prod_{j=1}^n \phi(z_j) \right\rangle_n^{(q,t)}
= \exp \( \sum_{k=1}^\infty kc(k)c(-k) \frac{1-q^k}{1-t^k} \).
$$
\end{thm}

\begin{proof}
First we see that
\begin{align*}
& \prod_{j=1}^n \phi(z_j)= e^{n c(0)} \prod_{k=1}^\infty \exp(c(k) p_k(\bz))
\exp(c(-k) \overline{p_k(\bz)}) \\
=& e^{n c(0)} \prod_{k=1}^\infty \( \sum_{a=0}^\infty \frac{c(k)^{a}}{a!} p_{(k^{a})}(\bz) \)
\( \sum_{b=0}^\infty \frac{c(-k)^{b}}{b!} \overline{p_{(k^{b})}(\bz)} \) \\
=& e^{n c(0)} \sum_{(1^{a_1} 2^{a_2} \cdots )} \sum_{(1^{b_1} 2^{b_2} \cdots )}
\( \prod_{k=1}^\infty \frac{c(k)^{a_k}c(-k)^{b_k}}{a_k! \, b_k!} \) 
p_{(1^{a_1}2^{a_2} \cdots )}(\bz)\overline{p_{(1^{b_1}2^{b_2} \cdots )}(\bz)},
\end{align*} 
where both $(1^{a_1}2^{a_2} \cdots )$ and $(1^{b_1}2^{b_2} \cdots )$
run over all partitions.
Therefore we have
$$
e^{-n c(0)} 
\left\langle \prod_{j=1}^n \phi(z_j) \right\rangle_n^{(q,t)} 
= \sum_{(1^{a_1} 2^{a_2} \cdots )} \sum_{(1^{b_1} 2^{b_2} \cdots )}
\( \prod_{k=1}^\infty \frac{c(k)^{a_k}}{a_k!} \frac{c(-k)^{b_k}}{b_k!} \) 
\frac{ \langle p_{(1^{a_1} 2^{a_2} \cdots )}, p_{(1^{b_1} 2^{b_2} \cdots )} 
\rangle_{\Delta^{\mr{Mac}}} }
{ \langle 1, 1 \rangle_{\Delta^{\mr{Mac}}} }.
$$
We recall the asymptotic behavior 
$$
\frac{ \langle p_{(1^{a_1} 2^{a_2} \cdots )}, p_{(1^{b_1} 2^{b_2} \cdots )} 
\rangle_{\Delta^{\mr{Mac}}} } 
{ \langle 1, 1 \rangle_{\Delta^{\mr{Mac}}} }  \qquad \longrightarrow  \qquad
\prod_{k=1}^\infty \delta_{a_k b_k} k^{a_k} a_k! \( \frac{1-q^k}{1-t^k}\)^{a_k}
$$
in the limit as $n \to \infty$, see \cite[Chapter VI (9.9) and (1.5)]{Mac}.
It follows from this that
\begin{align*}
&\lim_{n \to \infty} e^{-n c(0)} 
\left\langle \prod_{j=1}^n \phi(z_j) \right\rangle_n^{(q,t)}
= \sum_{(1^{a_1} 2^{a_2} \cdots )} 
\prod_{k=1}^\infty \frac{(k c(k) c(-k))^{a_k}}{a_k!} \(\frac{1-q^k}{1-t^k}\)^{a_k} \\
=& \prod_{k=1}^\infty \( \sum_{a=0}^\infty \frac{(k c(k) c(-k)\frac{1-q^k}{1-t^k})^{a}}{a!}\) 
= \exp \( \sum_{k=1}^\infty k c(k) c(-k) \frac{1-q^k}{1-t^k}\).
\end{align*}
Here $\sum_{k=1}^\infty k c(k) c(-k) \frac{1-q^k}{1-t^k}$ converges absolutely 
by the second assumption in \eqref{Eq:AssumptionSzego} and the Cauchy-Schwarz inequality,
because $| \frac{1-q^k}{1-t^k}| \le \frac{1+|q|^k}{1-|t|^{k}} \le 1+|q|$.
\end{proof}

Note that the present proof is similar to the corresponding one in \cite{BD}.
The result in \cite{BD} is the special case of Theorem \ref{Thm:SzegoMacdonald} with $q=t$.
As an example of this theorem,
the asymptotic behavior of the moment of $|\Psi^{\mr{A}}(\bz;\eta)|$
is given as follows.
A further asymptotic result is given by Corollary \ref{AsymMomentQ} below.

\begin{example} \label{ExampleMomentLimit}
Let $\gamma \in \bR$ and let $\eta$ be a complex number such that $|\eta| < 1$.
Then we have
$$
\lim_{n \to \infty} \left\langle |\Psi^{\mr{A}}(\bz;\eta)|^{2\gamma} \right\rangle_n^{(q,t)}
= \( \frac{(q |\eta|^2;t)_{\infty}}{(|\eta|^2;t)_{\infty}} \)^{\gamma^2}.
$$
This result is obtained by applying 
Theorem \ref{Thm:SzegoMacdonald} to $\phi(z)= |1+\eta z|^{2\gamma}$.
Then the Fourier coefficients of $\log \phi$ are
$c(k)=(-1)^{k-1} \eta^k \gamma/k$ and $c(-k)=(-1)^{k-1} \overline{\eta}^k \gamma/k$
for $k >0$, and $c(0)=0$. \qed
\end{example}

%
%%%%%%%%%%%%%%%%%%%%%%%%%%%%%%%%%%%%%%%%%%%%%%%%%%%%%%%%%%%%%%%%%%%%%%%%%%%%%%%%%%%%%%%%%%%%%%%%
\section{Circular ensembles and its $q$-analogue} \label{sectionCBEq}
%%%%%%%%%%%%%%%%%%%%%%%%%%%%%%%%%%%%%%%%%%%%%%%%%%%%%%%%%%%%%%%%%%%%%%%%%%%%%%%%%%%%%%%%%%%%%%%%
%

\subsection{Special case: $t=q^{\beta/2}$} 

In this subsection, we examine the results of the  last section for the special case
$t=q^{\beta/2}$ with $\beta>0$, i.e.,
we consider the weight function $\Delta^{\mr{Mac}}(\bz;q,q^{\beta/2})$.
Denote by $\langle \cdot \rangle_{n,\beta}^q$ the corresponding average.
Define the $q$-gamma function
(see e.g. \cite[(10.3.3)]{AAR}) by
$$
\Gamma_q(x)= (1-q)^{1-x} \frac{(q;q)_\infty}{(q^x;q)_{\infty}}.
$$

\begin{prop} \label{ThmMomentGamma}
Let $\beta$ be a positive real number.
For a positive integer $k$ and $\xi \in \bT$, we have
\begin{align} 
\left\langle \prod_{i=1}^k | \Psi(\bz;q^{i-1/2}\xi )|^2 \right\rangle_{n,\beta}^q
=& \prod_{i=0}^{k-1} \frac{\Gamma_t(\frac{2}{\beta}(i+1)) \Gamma_t(n+\frac{2}{\beta}(k+i+1))}
{\Gamma_t(\frac{2}{\beta}(k+i+1)) \Gamma_t(n+\frac{2}{\beta}(i+1))} \qquad
\text{(with $t=q^{\beta/2}$)} 
\label{eqCBEmomentQ} \\
=& \prod_{j=0}^{n-1} \frac{\Gamma_q (\frac{\beta}{2}j +2k+1) \Gamma_q(\frac{\beta}{2} j+1)}
{\Gamma_q(\frac{\beta}{2}j+k+1)^2}. \notag
\end{align}
\end{prop}

\begin{proof}
The claim follows immediately
from Corollary  \ref{CorMomentValue}
and the functional equation 
$\Gamma_q(1+x) = \frac{1-q^x}{1-q} \Gamma_q(x)$.
\end{proof}

Consider now the asymptotic behavior of this average in the limit as $n \to \infty$.
Put $[n]_q = (1-q^n)/(1-q)$.

\begin{cor} \label{AsymMomentQ}
For a positive integer $k$ and $\xi \in \bT$,
it holds that 
\begin{equation} \label{eqCBEmomentLimit}
\lim_{n \to \infty} ([n]_t)^{-2k^2/\beta}
\left\langle \prod_{i=1}^k | \Psi(\bz;q^{i-1/2}\xi )|^2 \right\rangle_{n,\beta}^q
= \prod_{i=0}^{k-1} \frac{\Gamma_t(\frac{2}{\beta}(i+1))}
{\Gamma_t(\frac{2}{\beta}(k+i+1))} \qquad \text{with $t=q^{\beta/2}$}.
\end{equation}
\end{cor}

\begin{proof}
Verify that 
\begin{equation} \label{eq:GammaQasym}
\lim_{n \to \infty} \frac{\Gamma_t(n+a)}{\Gamma_t(n) ([n]_t)^a} =1
\end{equation}
for any constant $a$.
Then the claim is clear from expression \eqref{eqCBEmomentQ}.
\end{proof}

\begin{example} \label{ExFq}
Denote by $\mcal{F}_\beta^q(k)$ the right-hand side of equation 
\eqref{eqCBEmomentLimit}.
Then we obtain
\begin{align}
\mcal{F}_{1}^q(k) =& \prod_{j=0}^{k-1} \frac{[2j+1]_{q^{\frac{1}{2}}} !}{[2k+2j+1]_{q^{\frac{1}{2}}}!},
\label{eqf1q} \\
\mcal{F}_{2}^q(k) =& \prod_{j=0}^{k-1} \frac{[j]_{q} !}{[j+k]_q!},  \label{eqf2q} \\
\mcal{F}_{4}^q(2k) =& 
\frac{([2]_q)^{2k^2}}{[2k-1]_q !!} \prod_{j=1}^{2k-1} \frac{[j]_q!}{[2j]_q!}.
\label{eqf4q}
\end{align}
Here $[n]_q!=[n]_q [n-1]_q \cdots [1]_q$ and $[2k-1]_q!! = [2k-1]_q [2k-3]_q \cdots [3]_q [1]_q$.
Equalities \eqref{eqf1q} and \eqref{eqf2q} are trivial because $\Gamma_q(n+1)=[n]_q!$.
We check relation \eqref{eqf4q}.
By definition, we have 
$$
\mcal{F}_4^q(2k)= 
\prod_{i=0}^{2k-1} \frac{\Gamma_{q^2}(\frac{1}{2}(i+1))}
{\Gamma_{q^2}(k+\frac{1}{2}(i+1))}
= \prod_{p=0}^{k-1} 
\frac{\Gamma_{q^2} (p+\frac{1}{2}) \Gamma_{q^2}(p+1)}
{\Gamma_{q^2} (k+p+\frac{1}{2}) \Gamma_{q^2}(k+p+1)}.
$$
Using the $q$-analogue of the Legendre duplication formula (see e.g. \cite[Theorem 10.3.5(a)]{AAR})
$$
\Gamma_q(2x) \Gamma_{q^2}(1/2) = (1+q)^{2x-1} \Gamma_{q^2}(x) \Gamma_{q^2}(x+1/2),
$$
we have
$$
\mcal{F}_4^q(2k)= \prod_{p=0}^{k-1} \frac{(1+q)^{2k} \Gamma_q(2p+1)}{\Gamma_q(2k+2p+1)}=
([2]_q)^{2k^2} \prod_{p=0}^{k-1} \frac{[2p]_q! }{[2k+2p]_q !}.
$$
Expression \eqref{eqf4q} can then be proven by  induction on $n$.
\end{example}

\subsection{Circular $\beta$-ensembles and Jack polynomials} \label{subsectionCBEJack}

We take the limit as $q \to 1$ of the results of the previous subsection.
Recall the formula
$$
\lim_{q \to 1} \frac{(q^a x;q)_{\infty}}{(x;q)_{\infty}} =(1-x)^{-a}
$$
for $|x|<1$ and $a \in \bR$, see \cite[Theorem 10.2.4]{AAR} for example.
Then we have  
$$
\lim_{q \to 1} \Delta^{\mr{Mac}}(\bz;q,q^{\beta/2}) 
= \prod_{1 \le i<j \le n} |z_i-z_j|^\beta =: \Delta^{\mr{Jack}}(\bz;2/\beta),
$$
which is a constant times the p.d.f. for Dyson's circular $\beta$-ensembles
(see \S 6). 
Denote by $\langle \cdot \rangle_{n,\beta}$ the corresponding average, i.e., 
for a function $f$ on $\bT^n$  define 
$$
\langle f \rangle_{n,\beta} = \lim_{q \to 1} \langle f \rangle_{n,\beta}^q
= \frac{\int_{\bT^n} f(\bz) \prod_{1 \le i<j \le n} |z_i-z_j|^\beta \dd \bz}
{\int_{\bT^n} \prod_{1 \le i<j \le n} |z_i-z_j|^\beta \dd \bz}.
$$

Let $\alpha >0$.
The Jack polynomial $P^{\mr{Jack}}_\lambda(x_1,\dots,x_n;\alpha)$ 
for each partition $\lambda$
is defined by the limit approached by the corresponding
  Macdonald polynomial,
$$
P^{\mr{Jack}}_\lambda(x_1,\dots,x_n;\alpha)
= \lim_{q \to 1} P^{\mr{Mac}}_\lambda(x_1,\dots,x_n;q,q^{1/\alpha}),
$$
see \cite[Chapter VI-10]{Mac} for detail.
Jack polynomials are orthogonal polynomials with respect to 
the weight function $\Delta^{\mr{Jack}}(\bz;\alpha)$.
In particular, $s_{\lambda}(x_1,\dots,x_n)=P^{\mr{Jack}}_\lambda(x_1,\dots,x_n;1)$
are called Schur polynomials,
and are irreducible characters of $U(n)$ associated with $\lambda$.

From the theorems in the last section, we have the following:
from Theorem \ref{ThmAverageMac}, we see that
\begin{equation} \label{AverageProductA}
\left\langle \prod_{l=1}^L \Psi^{\mr{A}}(\bz^{-1};\eta_l^{-1}) \cdot \prod_{k=1}^K
\Psi^{\mr{A}}(\bz;\eta_{L+k}) \right\rangle_{n,\beta} =
(\eta_1 \cdots \eta_L)^{-n}
\cdot P_{(n^L)}^{\mr{Jack}} (\eta_1, \dots, \eta_{L+K};\beta/2).
\end{equation}
For a positive real number $\gamma$ and complex number $\eta$ with $|\eta|<1$, 
we have from  Proposition \ref{PropMomentHypergeometric} that
\begin{equation} \label{MomentHypergeometricJack}
\left\langle |\Psi(\bz;\eta)|^{2\gamma} \right\rangle_{n, \CBE{}}
= {_2 F_1}^{(2/\beta)}(-\gamma, -\gamma; \frac{\beta}{2}(n-1)+1; 
|\eta|^2, \dots, |\eta|^2),
\end{equation}
where ${_2 F_1}^{(\alpha)}(a,b; c; x_1,\dots,x_n)$ is the hypergeometric function associated 
with Jack polynomials \cite{Kaneko1} defined by
$$
{_2 F_1}^{(\alpha)}(a,b; c; x_1,\dots,x_n)= \sum_{\lambda} 
\frac{[a]^{(\alpha)}_\lambda [b]^{(\alpha)}_\lambda}{[c]^{(\alpha)}_\lambda}
\frac{\alpha^{|\lambda|} P_\lambda^{\mr{Jack}}(x_1,\dots,x_n;\alpha)}{c'_\lambda(\alpha)}
$$
with 
$$
[u]_\lambda^{(\alpha)}=\prod_{s \in \lambda} (u-l'(s)/\alpha +a'(s)), \qquad 
\text{and} \qquad
c'_\lambda(\alpha)= \prod_{s \in \lambda}(\alpha(a(s)+1)+l(s)).
$$
For a positive integer $k$, and $\xi \in \bT$, by Theorem \ref{ThmMomentGamma} 
and Corollary \ref{AsymMomentQ} it holds that
\begin{equation} \label{MomentAsymptoticA}
\left\langle  | \Psi^{\mr{A}}(\bz;\xi )|^{2k} \right\rangle_{n,\beta}
= \prod_{i=0}^{k-1} \frac{\Gamma(\frac{2}{\beta}(i+1)) \Gamma(n+\frac{2}{\beta}(k+i+1))}
{\Gamma(\frac{2}{\beta}(k+i+1)) \Gamma(n+\frac{2}{\beta}(i+1))}
\sim 
\prod_{i=0}^{k-1} \frac{\Gamma(\frac{2}{\beta}(i+1)) }
{\Gamma(\frac{2}{\beta}(k+i+1))} \cdot n^{2k^2/\beta}
\end{equation}
in the limit as $n \to \infty$.
For a function
$\phi(z)=\exp(\sum_{k \in \bZ} c(k) z^k)$ on $\bT$ satisfying inequalities \eqref{Eq:AssumptionSzego},
 by Theorem \ref{Thm:SzegoMacdonald} it holds that
\begin{equation} \label{eq:SzegoJack}
\lim_{n \to \infty} e^{-n c(0)} 
\left\langle \prod_{j=1}^n \phi(z_j) \right\rangle_{n, \beta}
= \exp \( \frac{2}{\beta}\sum_{k=1}^\infty kc(k)c(-k) \). 
\end{equation}
In particular, for $\gamma \in \bR$  and a complex number  
$\eta$ such that $|\eta|< 1$,
we have
$$
\lim_{n \to \infty} \left\langle |\Psi^{\mr{A}}(\bz;\eta)|^{2\gamma} \right\rangle_{n,\beta}
= (1-|\eta|^2)^{-2 \gamma^2/\beta}.
$$

Several observations may be made concerning the above identities:  
equation \eqref{MomentHypergeometricJack} is
obtained by verifying the limits
$$
\lim_{t \to 1} \frac{(q^a)_\lambda^{(q,t)}}{(1-t)^{|\lambda|}} 
=\alpha^{|\lambda|} [a]_\lambda^{(\alpha)}, \qquad
\lim_{t \to 1} \frac{c'_\lambda(q,t)}{(1-t)^{|\lambda|}} =c'_\lambda(\alpha),
$$
with $q=t^\alpha$.
The expression for the moment is obtained in \cite{FK} using a different proof,
which employs a Selberg type integral evaluation.
Equation \eqref{MomentAsymptoticA} is also obtained in \cite{KS_zetafunctions} essentially
by the Selberg integral evaluation.
When $\beta=2$, equation \eqref{eq:SzegoJack} presents the strong Szeg\"{o} limit theorem for
a Toeplitz determinant.
Indeed, the average of the left-hand side of \eqref{eq:SzegoJack}
is then equal to the Toeplitz determinant $\det(d_{i-j})_{1 \le i,j \le n}$ of $\phi$,
where $d_i$ are Fourier coefficients of $\phi$.
Equation \eqref{eq:SzegoJack} with general $\beta>0$ is seen in \cite{Johansson1, Johansson2},
 but it may be noted that the present proof, employing symmetric function theory, is straightforward.
This expression is applied in \cite{Hyper} in order to observe an asymptotic behavior 
for Toeplitz  `hyperdeterminants'.

%
%%%%%%%%%%%%%%%%%%%%%%%%%%%%%%%%%%%%%%%%%%%%%%%%%%%%%%%%%%%%%%%%%%%%%%%%%%%%%%%%%%%%%%%%%%%%%%%%
\section{Jacobi polynomials due to Heckman and Opdam}
%%%%%%%%%%%%%%%%%%%%%%%%%%%%%%%%%%%%%%%%%%%%%%%%%%%%%%%%%%%%%%%%%%%%%%%%%%%%%%%%%%%%%%%%%%%%%%%%
%

The results obtained in \S\ref{sectionMacAverage} and \S\ref{sectionCBEq} will  be applied to random matrix 
polynomials from symmetric spaces of the type A root system in the next section.
In order to evaluate the corresponding polynomials of the BC type root system,
we here recall Heckman and Opdam's Jacobi polynomials and 
give some identities corresponding to \eqref{AverageProductA} and \eqref{MomentAsymptoticA}.

The dominance ordering associated with the root system of type BC is defined as follows:
for two partitions $\lambda=(\lambda_1,\lambda_2,\dots)$ and $\mu=(\mu_1,\mu_2,\dots)$,
$$
\mu \le \lambda \qquad \Leftrightarrow \qquad 
\mu_1 + \cdots+\mu_i \le \lambda_1+ \cdots +\lambda_i \quad \text{for all $i \ge 1$}.
$$
Let $\bC[\bx^{\pm 1}] = \bC[x_1^{\pm 1}, \dots, x_n^{\pm 1}]$ be the ring of all 
Laurent polynomials in $n$ variables $\bx=(x_1,\dots,x_n)$.
The Weyl group $W=\bZ_2 \wr \mf{S}_n = \bZ_2^n \rtimes \mf{S}_n$ of type $BC_n$
acts naturally on $\bZ^n$ and $\bC[\bx^{\pm 1}]$, respectively.
Denote by $\bC[\bx^{\pm 1}]^W$ the subring of all $W$-invariants in $\bC[\bx^{\pm 1}]$.
Let $\Delta^{\mr{HO}}(\bz;k_1,k_2,k_3)$ be a function on $\bT^n$ defined by 
$$
     \Delta^{\mr{HO}}(\bz;k_1,k_2,k_3) 
  =   \prod_{1 \le i < j \le n}
        |1-z_i z_j^{-1}|^{2 k_3} |1-z_i z_j|^{2 k_3} 
      \cdot 
     \prod_{1 \le j  \le n}
       |1-z_j|^{2k_1} |1-z_j^2|^{2k_2}.
$$
Here we assume $k_1$, $k_2$, and $k_3$ are real numbers such that
$$
     k_1+k_2>-1/2, \quad k_2 > -1/2, \quad k_3 \ge 0.
$$
Define an inner product on $\bC[\bx^{\pm 1}]^W$ by
$$
\langle f,g \rangle_{\Delta^{\mr{HO}}} = 
\frac{1}{2^n n!} \int_{\bT^n} f(\bz) g(\bz^{-1}) \Delta^{\mr{HO}}(\bz;k_1,k_2,k_3)
\dd \bz.
$$

For each partition $\mu$, we let 
$$
m^{\mr{BC}}_{\mu}(\bx)=\sum_{\nu \in W \mu} x_1^{\nu_1} \cdots x_n^{\nu_n},
$$
where $W\mu$ is the $W$-orbit of $\mu$
(cf. \eqref{eq:monomialA}).
These polynomials form a $\bC$-basis of $\bC[\bx^{\pm 1}]^W$.
Then, there exists a unique family of polynomials
$P^{\mr{HO}}_{\lambda}= P^{\mr{HO}}_{\lambda}(\bx;k_1,k_2,k_3) \in \bC[\bx^{\pm 1}]^W$
($\lambda$ are partitions such that $\ell(\lambda) \le n$)
satisfying two conditions: 
$$
P^{\mr{HO}}_{\lambda}(\bx)= m_{\lambda}^{\mr{BC}}(\bx)+ \sum_{\mu: \mu < \lambda} 
u_{\lambda \mu} m^{\mr{BC}}_{\mu}(\bx), 
\quad \text{with $u_{\lambda \mu} \in \bC$},
\qquad\qquad  \langle P^{\mr{HO}}_{\lambda}, P_{\mu}^{\mr{HO}} \rangle_{\Delta^{\mr{HO}}} = 0
\quad  
\text{if $\lambda \not= \mu$}.
$$
The Laurent polynomials $P_{\lambda}$ are known as Jacobi polynomials
associated with the root system of type $BC_n$ due to Heckman and Opdam,
see e.g. \cite{Diejen, Heckman, Mimachi}.
They can be seen as BC-analogues of Jack polynomials.

For a function $f$ on $\bT^n$, we denote by $\langle f \rangle_{n}^{k_1,k_2,k_3}$
the mean value of $f$ with respect to 
the weight function $\Delta^{\mr{HO}}(\bz;k_1,k_2,k_3)$:
$$
     \langle f \rangle_{n}^{k_1,k_2,k_3}
   =  \frac{\int_{\bT^n} f(\bz) \Delta^{\mr{HO}}(\bz;k_1,k_2,k_3)\dd \bz}{
  \int_{\bT^n} \Delta^{\mr{HO}}(\bz;k_1,k_2,k_3)\dd \bz}.  
$$
From the three parameters $k_1,k_2,k_3$, we define new parameters
$$
\tilde{k}_1 = k_1/k_3, \qquad \tilde{k}_2=(k_2+1)/k_3-1, \qquad \tilde{k}_3=1/k_3.
$$
Put
$$
\Psi^{\mr{BC}}(\bz;x)= \prod_{j=1}^n(1+x z_j)(1+x z_j^{-1}).
$$

\begin{thm} \label{Thm:MainTheorem}
The following relation holds 
\begin{equation} \label{eq:MainEq}
       \left\langle
          \Psi^{\mr{BC}}(\bz;x_1)\Psi^{\mr{BC}}(\bz;x_2) \cdots \Psi^{\mr{BC}}(\bz;x_m)
       \right\rangle_{n}^{k_1,k_2,k_3}     
     = (x_1 \cdots x_m)^n
       P^{\mr{HO}}_{(n^m)}(x_1,\dots,x_m;\tilde{k}_1,\tilde{k}_2,\tilde{k}_3).
\end{equation}
\end{thm}

In order to prove this,
we need the following dual Cauchy identity obtained by Mimachi \cite{Mimachi}.

\begin{prop}[\cite{Mimachi}] \label{Thm:Mimachi}
Let $\bx=(x_1,\dots, x_n)$ and 
$\by=(y_1,\dots, y_m)$ be sequences of indeterminates.
Jacobi polynomials $P^{\mr{HO}}_{\lambda}$ satisfy the equality 
$$
\prod_{i=1}^n \prod_{j=1}^m (x_i+x_i^{-1} - y_j - y_j^{-1}) 
= \sum_{\lambda \subset (m^n)} (-1)^{|\tilde{\lambda}|} 
P^{\mr{HO}}_{\lambda}(\bx;k_1,k_2,k_3) 
P^{\mr{HO}}_{\tilde{\lambda}}(\by;\tilde{k}_1,\tilde{k}_2,\tilde{k}_3),
$$
where $\tilde{\lambda}=(n-\lambda_m', n-\lambda_{m-1}', \dots, n-\lambda_1')$.
\qed
\end{prop}

\begin{proof}[Proof of Theorem \ref{Thm:MainTheorem}]
We see that
$$
       \Psi^{\mr{BC}}(\bz;x_1) \Psi^{\mr{BC}}(\bz;x_2) \cdots \Psi^{\mr{BC}}(\bz;x_m) 
    =  (x_1 \cdots x_m)^n \prod_{i=1}^m 
        \prod_{j=1}^n (x_i + x_i^{-1} + z_j + z_j^{-1}).
$$
Using Proposition \ref{Thm:Mimachi} we have
\begin{align*}
     &  \left\langle
          \Psi^{\mr{BC}}(\bz;x_1)\Psi^{\mr{BC}}(\bz;x_2) \cdots \Psi^{\mr{BC}}(\bz;x_m)
       \right\rangle_{n}^{k_1,k_2,k_3} \\
   = & (x_1 \cdots x_m)^n \sum_{\lambda \subset (m^n)} 
        P^{\mr{HO}}_{\tilde{\lambda}}(x_1,\dots, x_m;\tilde{k}_1,\tilde{k}_2,\tilde{k}_3) 
        \langle P^{\mr{HO}}_{\lambda}(\bz;k_1,k_2,k_3) \rangle_{n}^{k_1,k_2,k_3}.
\end{align*}
By the orthogonality relation for Jacobi polynomials, 
we have
$$
        \langle P^{\mr{HO}}_{\lambda}(\bz;k_1,k_2,k_3) \rangle_{n}^{k_1,k_2,k_3}
    =   \begin{cases}  
          1, & \text{if $\lambda = (0)$}, \\ 0, & \text{otherwise},
        \end{cases}
$$
and we thus obtain the theorem.
\end{proof}

\begin{remark}
Using Theorem 2.1 in \cite{Mimachi},
we derive a more general form of equation \eqref{eq:MainEq}
including a Macdonald-Koornwinder polynomial.
\end{remark}

\begin{thm} \label{Thm:Main2}
Let
$$
\mcal{F}(m;k_1,k_2,k_3)= \prod_{j=0}^{m-1} \frac{\sqrt{\pi}}{2^{k_1 +2 k_2+j k_3-1}
\Gamma(k_1+k_2+\frac{1}{2}+j k_3)}.
$$
The $m$-th moment of $\Psi^{\mr{BC}}(\bz;1)$ is given by 
$$
       \left\langle
          \Psi^{\mr{BC}}(\bz;1)^m
       \right\rangle_{n}^{k_1,k_2,k_3}  
   = \mcal{F}(m;\tilde{k}_1,\tilde{k}_2,\tilde{k}_3) \cdot \prod_{j=0}^{m-1} 
\frac{\Gamma(n+ \tilde{k}_1+2\tilde{k}_2+j \tilde{k}_3 ) 
\Gamma(n+ \tilde{k}_1+\tilde{k}_2+\frac{1}{2}+j \tilde{k}_3 )}
{\Gamma(n+ \frac{\tilde{k}_1}{2}+\tilde{k}_2+\frac{j \tilde{k}_3}{2} )
\Gamma(n+ \frac{\tilde{k}_1}{2}+\tilde{k}_2+\frac{1+j \tilde{k}_3}{2} )}. 
$$
\end{thm}

\begin{proof}
By Theorem \ref{Thm:MainTheorem} we have
\begin{equation} \label{eq:MTspecial}
   \left\langle
          \Psi^{\mr{BC}}(\bz;1)^m
       \right\rangle_{n}^{k_1,k_2,k_3}
= P^{\mr{HO}}_{(n^m)}(1^m;\tilde{k}_1,\tilde{k}_2,\tilde{k}_3).
\end{equation}
The special case $P_{\lambda}^{\mr{HO}}(1,1,\dots,1;k_1,k_2,k_3)$ is known
and is given as follows (see e.g. \cite{Diejen}
\footnote{The connection between ours notation and van Diejen's \cite{Diejen} is given by
$\nu_0 = k_1+k_2, \ \nu_1=k_2, \ \nu=k_3$.}):
for a partition $\lambda$ of length $\le m$,
\begin{align*}
P^{\mr{HO}}_{\lambda}(\underbrace{1, \dots, 1}_m;k_1,k_2,k_3)
=&  2^{2|\lambda|} \prod_{1 \le i <j \le m}
   \frac{(\rho_i+ \rho_j+k_3)_{\lambda_i+\lambda_j} 
          (\rho_i- \rho_j+k_3)_{\lambda_i-\lambda_j}}
        {(\rho_i+ \rho_j)_{\lambda_i+\lambda_j} 
          (\rho_i- \rho_j)_{\lambda_i-\lambda_j}}  \\
  & \quad \times \prod_{j=1}^m 
     \frac{(\frac{k_1}{2} +k_2 + \rho_j)_{\lambda_j} (\frac{k_1+1}{2} + \rho_j)_{\lambda_j}}
          {(2 \rho_j)_{2 \lambda_j}} 
\end{align*}
with
$\rho_j= (m-j)k_3 + \frac{k_1}{2}+k_2$.
Here $(a)_n = \Gamma(a+n) / \Gamma(a)$ is the Pochhammer symbol.
Substituting $(n^m)$ for $\lambda$, we have 
\begin{align}
& P^{\mr{HO}}_{(n^m)} (1^m; k_1,k_2,k_3)  \notag \\
   =& \prod_{1 \le i <j \le m} \frac{(k_1+2 k_2+(2m-i-j+1)k_3)_{2n}}
                                    {(k_1+2 k_2+(2m-i-j)k_3)_{2n}}
    \cdot \prod_{j=0}^{m-1} 
      \frac{2^{2n} (k_1+2 k_2+j k_3)_n (k_1+k_2+\frac{1}{2}+j k_3)_n}
           {(k_1 +2 k_2+2 j k_3)_{2n}}. \label{eq:moment_product1}
\end{align}

A simple algebraic manipulation of the first product on the right-hand side of   \eqref{eq:moment_product1} yields   
$$
\prod_{1 \le i <j \le m} \frac{(k_1+2 k_2+(2m-i-j+1)k_3)_{2n}}
                                    {(k_1+2 k_2+(2m-i-j)k_3)_{2n}}
= \prod_{j=0}^{m-1} \frac{(k_1+2k_2+ 2jk_3)_{2n}}{(k_1+2k_2+j k_3)_{2n}}
$$
and therefore we obtain
$$
P_{(n^m)}^{\mr{HO}} (1^m; k_1,k_2,k_3) =
\prod_{j=0}^{m-1} \frac{2^{2n} (k_1+k_2+\frac{1}{2}+j k_3)_n}{(n+k_1+2k_2+jk_3)_{n}}.
$$
Combining the above result with equation \eqref{eq:MTspecial}, we have 
\begin{equation} \label{eq:Main2}
       \left\langle
          \Psi^{\mr{BC}}(\bz;1)^m
       \right\rangle_{n}^{k_1,k_2,k_3}  
   = \prod_{j=0}^{m-1} 
\frac{2^{2n} \Gamma(n+ \tilde{k}_1+2\tilde{k}_2+j \tilde{k}_3 ) 
\Gamma(n+ \tilde{k}_1+\tilde{k}_2+\frac{1}{2}+j \tilde{k}_3 )}
{\Gamma( \tilde{k}_1+\tilde{k}_2+\frac{1}{2} +j \tilde{k}_3)
\Gamma(2n+ \tilde{k}_1+2\tilde{k}_2+j \tilde{k}_3 )}. 
\end{equation}

Finally, we apply the formula
$$
\Gamma(2a) = \frac{2^{2a-1}}{\sqrt{\pi}} \Gamma(a) \Gamma(a+\frac{1}{2})
$$
to $\Gamma(2n+\tilde{k}_1+2\tilde{k}_2+j \tilde{k}_3)$ in equation \eqref{eq:Main2}
and we then have the theorem.
\end{proof}

\begin{cor} \label{cor:Main}
It holds that
$$
      \left\langle
          \Psi^{\mr{BC}}(\bz;1)^m
       \right\rangle_{n}^{k_1,k_2,k_3}  
   \sim
 \mcal{F}(m;\tilde{k}_1,\tilde{k}_2,\tilde{k}_3) \cdot 
n^{m(\tilde{k}_1+\tilde{k}_2)+\frac{1}{2}m(m-1)\tilde{k}_3},
$$
as $n \to \infty$.
\end{cor}

\begin{proof}
The claim follows from the previous theorem and
the asymptotics of the gamma function (cf \eqref{eq:GammaQasym}):
$\Gamma(n+a) \sim \Gamma(n) n^a$ for a constant $a$.
\end{proof}

%
%%%%%%%%%%%%%%%%%%%%%%%%%%%%%%%%%%%%%%%%%%%%%%%%%%%%%%%%%%%%%%%%%%%%%%%%%%%%%%%%%%%%%%%%%%%%%%%%
\section{Random matrix ensembles associated with compact symmetric spaces}
%%%%%%%%%%%%%%%%%%%%%%%%%%%%%%%%%%%%%%%%%%%%%%%%%%%%%%%%%%%%%%%%%%%%%%%%%%%%%%%%%%%%%%%%%%%%%%%%
%

Finally, we apply the theorems obtained above to 
compact symmetric spaces as classified by Cartan.
These symmetric spaces are labeled  A I, BD I, C II, and so on, 
see e.g. Table 1 in \cite{CM}.
Let $G/K$ be such a compact symmetric space.
Here $G$ is a compact subgroup of $GL(N,\bC)$ for some positive integer $N$,
and $K$ is a closed subgroup of $G$.
Then the space $G/K$ is realized as the subset $S$ of $G$: $S \simeq G/K$
and the probability measure $\dd M$ on $S$ is induced from the quotient space $G/K$.
We consider $S$ as a probability space with the measure $\dd M$ and 
call the random matrix ensemble associated with $G/K$.
See \cite{Duenez} for details.

The random matrix ensembles considered in
\S\ref{subsectionA}, \S\ref{subsectionAI}, and \S\ref{subsectionAII}
are called Dyson's circular $\beta$-ensembles, see \cite{Dyson, Mehta}.
The identities in these subsections 
follow immediately from expressions \eqref{AverageProductA} and \eqref{MomentAsymptoticA} 
(see also Example \ref{ExFq}) .
Similarly, identities after \S\ref{subsectionB} follows from 
Theorem \ref{Thm:MainTheorem},  Theorem \ref{Thm:Main2}, and Corollary \ref{cor:Main}.

Note that the results in \S\ref{subsectionA}, \S\ref{subsectionB}, \S\ref{subsectionC}, 
and  \S\ref{subsectionD} are results for 
compact Lie groups (which are not proper symmetric spaces)
previously presented in \cite{BG}.

\subsection{$U(n)$ -- type A} \label{subsectionA}

Consider the unitary group $U(n)$ with the normalized Haar measure.
This space has a simple root system of type A.
The corresponding p.d.f. for eigenvalues $z_1,\dots,z_n$ of $M \in U(n)$
is proportional to $\Delta^{\mr{Jack}}(\bz;1)$.
This random matrix ensemble is called the circular unitary ensemble (CUE).

For complex numbers $\eta_1,\dots,\eta_L, \eta_{L+1},\dots, \eta_{L+K}$, 
it follows from equation \eqref{AverageProductA} that
\begin{align*}
& \left\langle \prod_{i=1}^L  \det(I+\eta_i^{-1} M^{-1}) \cdot \prod_{i=1}^K 
\det(I+\eta_{L+i} M) \right\rangle_{U(n)} \\
=&  \left\langle \prod_{i=1}^L  \Psi^{\mr{A}}(\bz^{-1};\eta_i^{-1}) \cdot \prod_{i=1}^K 
\Psi^{\mr{A}}(\bz;\eta_{L+i}) \right\rangle_{n,2} =
\prod_{i=1}^L \eta_i^{-n} \cdot s_{(n^L)} (\eta_1,\dots,\eta_{L+K}).
\end{align*}
In addition, from equation \eqref{MomentAsymptoticA} we obtain
$$
\left\langle |\det(I+ \xi M)|^{2m} \right\rangle_{U(n)}
=  \prod_{j=0}^{m-1}\frac{j! (n+j+m)!}{(j+m)! (n+j)!}
\sim
\prod_{j=0}^{m-1}\frac{j!}{(j+m)!} \cdot n^{m^2}
$$
for any $\xi \in \bT$.

\subsection{$U(n)/O(n)$ -- type A I} \label{subsectionAI}

Consider the ensemble $S(n)$ associated with the symmetric space $U(n)/O(n)$.
The space $S(n)$ is the set of all symmetric matrices in $U(n)$.
The corresponding p.d.f. for eigenvalues $z_1,\dots,z_n$
is proportional to $\Delta^{\mr{Jack}}(\bz;2) = \prod_{1 \le i<j \le n} |z_i-z_j|$.
This random matrix ensemble is called the circular orthogonal ensemble (COE).
We have
\begin{align*}
&\left\langle \prod_{i=1}^L  \det(I+\eta_i^{-1} M^{-1}) \cdot \prod_{i=1}^K 
\det(I+\eta_{L+i} M) \right\rangle_{S(n)} \\
=&  \left\langle \prod_{i=1}^L  \Psi^{\mr{A}}(\bz^{-1};\eta_i^{-1}) \cdot \prod_{i=1}^K 
\Psi^{\mr{A}}(\bz;\eta_{L+i}) \right\rangle_{n,1}
= \prod_{i=1}^L \eta_i^{-n} \cdot P_{(n^L)}^{\mr{Jack}} (\eta_1,\dots,\eta_{L+K};1/2).
\end{align*}
For $\xi \in \bT$,  we obtain
$$
\left\langle |\det(I+ \xi M)|^{2m} \right\rangle_{S(n)}
= \prod_{j=0}^{m-1} \frac{(2j+1)! (n+2m+2j+1)!}{(2m+2j+1)! (n+2j+1)!}
\sim
\prod_{j=0}^{m-1}\frac{(2j+1)!}{(2m+2j+1)!} \cdot n^{2m^2}.
$$

\subsection{$U(2n)/Sp(2n)$ -- type A II} \label{subsectionAII}

Consider the ensemble $S(n)$ associated with the symmetric space $U(2n)/Sp(2n)$.
The space $S(n)$ is the set of all self-dual matrices in $U(2n)$,
i.e., $M \in S(n)$ is a unitary matrix satisfying
$M=J \trans{M} \trans{J}$ with 
$J=\(\begin{smallmatrix} 0 & I_n \\ -I_n & 0 \end{smallmatrix} \)$.  
This random matrix ensemble is called the circular symplectic ensemble (CSE).
The eigenvalues of $M \in S(n)$  are of the form
$z_1,z_1,z_2,z_2,\dots,z_n,z_n$
and so the characteristic polynomial is given as
$\det(I+xM)= \prod_{j=1}^n (1+x z_j)^2$.
The corresponding p.d.f. for $z_1,\dots,z_n$
is proportional to $\Delta^{\mr{Jack}}(\bz;1/2)=\prod_{1 \le i<j \le n} |z_i-z_j|^4$.
We have
\begin{align*}
&\left\langle \prod_{i=1}^L  \det(I+\eta_i^{-1} M^{-1})^{1/2} \cdot \prod_{i=1}^K 
\det(I+\eta_{L+i} M)^{1/2} \right\rangle_{S(n)} \\
=& \left\langle \prod_{i=1}^L  \Psi^{\mr{A}}(\bz^{-1};\eta_i^{-1}) \cdot \prod_{i=1}^K 
\Psi^{\mr{A}}(\bz;\eta_{L+i}) \right\rangle_{n,4} 
= \prod_{i=1}^L x_i^{-n} \cdot P_{(n^L)}^{\mr{Jack}} (\eta_1,\dots,\eta_{L+K};2).
\end{align*}
For $\xi \in \bT$, we obtain
$$
\left\langle |\det(I+ \xi M)|^{2m} \right\rangle_{S(n)}
= \prod_{j=0}^{2m-1} \frac{\Gamma(\frac{j+1}{2}) \Gamma(n+m+\frac{j+1}{2})}
{\Gamma(m+\frac{j+1}{2}) \Gamma(n+\frac{j+1}{2})}
\sim
\frac{2^m}{(2m-1)!! \ \prod_{j=1}^{2m-1}(2j-1)!!} \cdot n^{2m^2}.
$$

\subsection{$SO(2n+1)$ -- type B} \label{subsectionB}

Consider the special orthogonal group $SO(2n+1)$.
An element $M$ in $SO(2n+1)$ is an  orthogonal matrix in $SL(2n+1,\bR)$, 
with eigenvalues given by
$z_1,z_1^{-1},\cdots, z_n,z_n^{-1},1$. 
From Weyl's integral formula,
the corresponding p.d.f.  of $z_1,z_2,\dots,z_n$ is 
proportional to
$\Delta^{\mr{HO}}(\bz;1,0,1)$,
and therefore it follows from Theorem \ref{Thm:MainTheorem} that
$$
\left\langle \prod_{i=1}^m \det(I+x_i M) \right\rangle_{SO(2n+1)}
= \prod_{i=1}^m (1+x_i) \cdot 
\left\langle \prod_{i=1}^m \Psi^{\mr{BC}}(\bz;x_i) \right\rangle_{n}^{1,0,1}
= \prod_{i=1}^m x_i^n (1+x_i) \cdot 
P_{(n^m)}^{\mr{HO}}(x_1,\dots,x_m;1,0,1).
$$ 
Here $P_\lambda^{\mr{HO}}(x_1,\dots,x_m;1,0,1)$ is just the irreducible character of $SO(2m+1)$
associated with the partition $\lambda$.
Theorem \ref{Thm:Main2}, Corollary \ref{cor:Main}, and a simple calculation lead to 
$$
\left\langle \det(I+ M)^m \right\rangle_{SO(2n+1)}
= 2^{m} \prod_{j=0}^{m-1} 
\frac{ \Gamma(2n+ 2j+2 ) }
{2^{j} (2j+1)!!  \ \Gamma(2n+ j+1)}
\sim \frac{2^{2m}}{\prod_{j=1}^{m} (2j-1)!!} n^{m^2/2+m/2}
$$
in the limit as $n \to \infty$.

\subsection{$Sp(2n)$ -- type C} \label{subsectionC}

Consider the symplectic group $Sp(2n)$,
i.e., a matrix $M \in Sp(2n)$ belongs to $U(2n)$ and satisfies 
$M J \trans{M}=J$, where
$J=\(\begin{smallmatrix} O_n & I_n \\ -I_n & O_n \end{smallmatrix} \)$.
The eigenvalues  are given by
$z_1,z_1^{-1},\cdots, z_n,z_n^{-1}$. 
The corresponding p.d.f.  of $z_1,z_2,\dots,z_n$ is 
proportional to
$\Delta^{\mr{HO}}(\bz;0,1,1)$
and therefore we have
$$
\left\langle \prod_{i=1}^m \det(I+x_i M) \right\rangle_{Sp(2n)}
= 
\left\langle \prod_{i=1}^m \Psi^{\mr{BC}}(\bz;x_i) \right\rangle_{n}^{0,1,1}
= \prod_{i=1}^m x_i^n \cdot 
P_{(n^m)}^{\mr{HO}}(x_1,\dots,x_m;0,1,1).
$$ 
Here $P_\lambda^{\mr{HO}}(x_1,\dots,x_m;0,1,1)$ is just the irreducible character of $Sp(2m)$
associated with the partition $\lambda$.
We obtain
$$
\left\langle \det(I+ M)^m \right\rangle_{Sp(2n)}
=  \prod_{j=0}^{m-1}
\frac{\Gamma(2n+2j+3) }{2^{j+1} \cdot  (2j+1)!!
 \ \Gamma(2n+j+2)} 
\sim  \frac{1}{ \prod_{j=1}^{m} (2j-1)!!} \cdot n^{m^2/2+m/2}.
$$

\subsection{$SO(2n)$ -- type D} \label{subsectionD}

Consider the special orthogonal group $SO(2n)$.
The eigenvalues of a matrix $M \in SO(2n)$ are of the form
$z_1,z_1^{-1},\cdots, z_n,z_n^{-1}$. 
The corresponding p.d.f.  of $z_1,z_2,\dots,z_n$ is 
proportional to
$\Delta^{\mr{HO}}(\bz;0,0,1)$,
and therefore we have
$$
\left\langle \prod_{i=1}^m \det(I+x_i M) \right\rangle_{SO(2n)}
= 
\left\langle \prod_{i=1}^m \Psi^{\mr{BC}}(\bz;x_i) \right\rangle_{n}^{0,0,1}
= \prod_{i=1}^m x_i^n \cdot 
P_{(n^m)}^{\mr{HO}}(x_1,\dots,x_m;0,0,1).
$$ 
Here $P_\lambda^{\mr{HO}}(x_1,\dots,x_m;0,0,1)$ is just the irreducible character of $O(2m)$
(not $SO(2m)$)
associated with the partition $\lambda$.
We  have
$$
\left\langle \det(I+ M)^m \right\rangle_{SO(2n)}
=  \prod_{j=0}^{m-1} \frac{\Gamma(2n+2j)}{2^{j-1} \, (2j-1)!! \ \Gamma(2n+j)} 
\sim  \frac{2^m}{\prod_{j=1}^{m-1} (2j-1)!!} \cdot n^{m^2/2-m/2}.
$$

\subsection{$U(2n+r)/(U(n+r)\times U(n))$ -- type A III} 

Let $r$ be a non-negative integer.
Consider the random matrix ensemble $G(n,r)$ associated with $U(2n+r)/(U(n+r)\times U(n))$.
The explicit expression of a matrix in $G(n,r)$ is omitted here, but may be found in \cite{Duenez}.
The eigenvalues of a matrix $M \in G(n,r) \subset U(2n+r)$ are  of the form
\begin{equation} \label{eq:Eigenvalues}
z_1,z_1^{-1},\cdots, z_n,z_n^{-1},\underbrace{1,1,\dots, 1}_r.
\end{equation}
The corresponding p.d.f.  of $z_1,z_2,\dots,z_n$ is proportional to
$\Delta^{\mr{HO}}(\bz;r,\frac{1}{2},1)$,
and therefore we have
$$
\left\langle \prod_{i=1}^m \det(I+x_i M) \right\rangle_{G(n,r)}
= \prod_{i=1}^m (1+x_i)^r \cdot 
\left\langle \prod_{i=1}^m \Psi^{\mr{BC}}(\bz;x_i) \right\rangle_{n}^{r,\frac{1}{2},1} 
= \prod_{i=1}^m (1+x_i)^r x_i^n \cdot 
P^{\mr{HO}}_{(n^m)}(x_1,\dots,x_m;r,\frac{1}{2},1).
$$
We  obtain
\begin{align*}
& \left\langle \det(I+ M)^m \right\rangle_{G(n,r)}
= 2^{mr} \left\langle \Psi^{\mr{BC}}(\bz;1) \right\rangle_{n}^{r,\frac{1}{2},1} \\
=&  \frac{\pi^{m/2}}{\prod_{j=0}^{m-1} 2^{j} (r+j)! }
\prod_{j=0}^{m-1} \frac{\Gamma(n+r+j+1)^2}
{\Gamma(n+\frac{r+j+1}{2}) \Gamma(n+\frac{r+j}{2}+1)} 
\sim  \frac{\pi^{m/2}}{2^{m(m-1)/2} \prod_{j=0}^{m-1} (r+j)! }
 \cdot n^{m^2/2 + rm}.
\end{align*}

\subsection{$O(2n+r)/(O(n+r) \times O(n))$ -- type BD I} 

Let $r$ be a non-negative integer.
Consider the random matrix ensemble $G(n,r)$ associated with the compact
symmetric space $O(2n+r)/(O(n+r) \times O(n))$.
The eigenvalues of a matrix $M \in G(n,r) \subset O(2n+r)$ are  of the form \eqref{eq:Eigenvalues}.
The corresponding p.d.f.  of $z_1,z_2,\dots,z_n$ is proportional to
$\Delta^{\mr{HO}}(\bz;\frac{r}{2},0,\frac{1}{2})$,
and therefore we have
$$
\left\langle \prod_{i=1}^m \det(I+x_i M) \right\rangle_{G(n,r)}
= \prod_{i=1}^m (1+x_i)^r \cdot 
\left\langle \prod_{i=1}^m \Psi^{\mr{BC}}(\bz;x_i) \right\rangle_{n}^{\frac{r}{2},0,\frac{1}{2}}
= \prod_{i=1}^m (1+x_i)^r x_i^n \cdot 
P_{(n^m)}^{\mr{HO}}(x_1,\dots,x_m;r,1,2).
$$ 
We  obtain
$$
\left\langle \det(I+ M)^m \right\rangle_{G(n,r)}
=  2^{mr}
\prod_{j=0}^{m-1} \frac{\Gamma(2n+4j+2r+3)}{2^{2j+r+1}(4j+2r+1)!! \
\Gamma(2n+2j+r+2)}
\sim \frac{2^{mr}}{\prod_{j=0}^{m-1}(4j+2r+1)!!} \cdot n^{m^2+rm}.
$$

\subsection{$Sp(2n)/U(n)$ -- type C I}

Consider the random matrix ensemble $S(n)$ associated with the compact
symmetric space $Sp(2n) /(Sp(2n)\cap SO(2n)) \simeq Sp(2n)/U(n)$.
The eigenvalues of a matrix $M \in S(n) \subset Sp(2n)$ are  of the form
$z_1,z_1^{-1},\cdots, z_n,z_n^{-1}$. 
The corresponding p.d.f.  of $z_1,z_2,\dots,z_n$ is proportional to
$\Delta^{\mr{HO}}(\bz;0,\frac{1}{2},\frac{1}{2})$,
and therefore we have
$$
\left\langle \prod_{i=1}^m \det(I+x_i M) \right\rangle_{S(n)}
=
\left\langle \prod_{i=1}^m \Psi^{\mr{BC}}(\bz;x_i) \right\rangle_{n}^{0,\frac{1}{2},\frac{1}{2}}
= \prod_{i=1}^m  x_i^n \cdot 
P_{(n^m)}^{\mr{HO}}(x_1,\dots,x_m;0,2,2).
$$ 
We  obtain
$$
\left\langle \det(I+ M)^m \right\rangle_{S(n)}
=  
\prod_{j=0}^{m-1} \frac{(n+2j+3) \Gamma(2n+4j+5)}{2^{2j+2}(4j+3)!! \
\Gamma(2n+2j+4)}
\sim \frac{1}{2^m \prod_{j=1}^m (4j-1)!!} \cdot n^{m^2+m}.
$$

\subsection{$Sp(4n+2r)/(Sp(2n+2r) \times Sp(2n))$ -- type C II} 

Let $r$ be a non-negative integer.
Consider the random matrix ensemble $G(n,r)$ associated with the compact
symmetric space $Sp(4n+2r)/(Sp(2n+2r) \times Sp(2n))$.
The eigenvalues of a matrix $M \in G(n,r) \subset Sp(4n+2r)$ are  of the form
$$
z_1,z_1,z_1^{-1},z_1^{-1},\cdots, z_n,z_n, z_n^{-1},z_n^{-1}, \underbrace{1,\dots,1}_{2r}.
$$ 
The corresponding p.d.f.  of $z_1,z_2,\dots,z_n$ is proportional to
$\Delta^{\mr{HO}}(\bz;2r,\frac{3}{2},2)$,
and therefore we have
\begin{align*}
\left\langle \prod_{i=1}^m \det(I+x_i M)^{1/2} \right\rangle_{G(n,r)}
=&\prod_{i=1}^m (1+x_i)^r
\left\langle \prod_{i=1}^m \Psi^{\mr{BC}}(\bz;x_i) \right\rangle_{n}^{2r,\frac{3}{2},2} \\
=& \prod_{i=1}^m  (1+x_i)^rx_i^n \cdot 
P^{\mr{HO}}_{(n^m)}(x_1,\dots,x_m;r,\frac{1}{4},\frac{1}{2}).
\end{align*}
We  obtain
\begin{align*}
\left\langle \det(I+ M)^m \right\rangle_{G(n,r)}
=&  
\frac{2^{4mr+m^2+m}}{\prod_{j=0}^{m-1} (4j+4r+1)!!} \cdot
\frac{\prod_{p=1}^{4m} \Gamma(n+r+\frac{p+1}{4})}{\prod_{j=1}^{2m}
\Gamma(n+\frac{r}{2}+\frac{j}{4}) \Gamma(n+\frac{r+1}{2}+\frac{j}{4})} \\
\sim & 
\frac{2^{4mr+m^2+m}}{\prod_{j=0}^{m-1} (4j+4r+1)!!} n^{m^2+2mr}.
\end{align*}

\subsection{$SO(4n+2)/U(2n+1)$ -- type D III-odd}

Consider the random matrix ensemble $S(n)$ associated with the compact
symmetric space $SO(4n+2)/(SO(4n+2) \cap Sp(4n+2)) \simeq  SO(4n+2)/U(2n+1)$.
The eigenvalues of a matrix $M \in S(n) \subset SO(4n+2)$ are  of the form
$z_1,z_1,z_1^{-1},z_1^{-1},\cdots, z_n,z_n, z_n^{-1},z_n^{-1}, 1,1$. 
The corresponding p.d.f.  of $z_1,z_2,\dots,z_n$ is proportional to
$\Delta^{\mr{HO}}(\bz;2,\frac{1}{2},2)$
and therefore we have
\begin{align*}
\left\langle \prod_{i=1}^m \det(I+x_i M)^{1/2} \right\rangle_{S(n)}
=&\prod_{i=1}^m (1+x_i)
\left\langle \prod_{i=1}^m \Psi^{\mr{BC}}(\bz;x_i) \right\rangle_{n}^{2,\frac{1}{2},2} \\
=& \prod_{i=1}^m  (1+x_i) x_i^n \cdot 
P^{\mr{HO}}_{(n^m)}(x_1,\dots,x_m;1,-\frac{1}{4},\frac{1}{2}).
\end{align*}
We  obtain
$$
\left\langle \det(I+ M)^m \right\rangle_{S(n)}
=  
\frac{2^{m^2+5m}}{\prod_{j=1}^{m} (4j-1)!!} \cdot
\prod_{j=1}^{2m} \frac{\Gamma(n+\frac{j}{2}+\frac{3}{4}) \Gamma(n+\frac{j}{2})}
{\Gamma(n+\frac{j}{4}) \Gamma(n+\frac{j}{4} +\frac{1}{2})}
\sim 
\frac{2^{m^2+5m}}{\prod_{j=1}^{m} (4j-1)!!} \cdot n^{m^2+m}.
$$

\subsection{$SO(4n)/U(2n)$ -- type D III-even}

Consider the random matrix ensembles $S(n)$ associated  with the compact
symmetric space $SO(4n)/(SO(4n) \cap Sp(4n)) \simeq SO(4n)/U(2n)$.
The eigenvalues of the matrix $M \in S(n) \subset SO(4n)$ are  of the form
$$
z_1,z_1,z_1^{-1},z_1^{-1},\cdots, z_n,z_n, z_n^{-1},z_n^{-1}.
$$ 
The corresponding p.d.f.  of $z_1,z_2,\dots,z_n$ is proportional to
$\Delta^{\mr{HO}}(\bz;0,\frac{1}{2},2)$
and therefore we have
$$
\left\langle \prod_{i=1}^m \det(I+x_i M)^{1/2} \right\rangle_{S(n)}
=
\left\langle \prod_{i=1}^m \Psi^{\mr{BC}}(\bz;x_i) \right\rangle_{n}^{0,\frac{1}{2},2} 
=  P_{(n^m)}^{\mr{HO}}(x_1,\dots,x_m;0,-\frac{1}{4},\frac{1}{2}).
$$
Hence we  obtain
$$
\left\langle \det(I+ M)^m \right\rangle_{S(n)}
=  
\frac{2^{m^2+m}}{\prod_{j=1}^{m-1} (4j-1)!!} \cdot
\prod_{j=0}^{2m-1} \frac{\Gamma(n+\frac{j}{2}+\frac{1}{4}) \Gamma(n+\frac{j-1}{2})}
{\Gamma(n+\frac{j-1}{4}) \Gamma(n+\frac{j+1}{4})}
\sim 
\frac{2^{m^2+m}}{\prod_{j=1}^{m-1} (4j-1)!!} \cdot n^{m^2-m}.
$$

%
%%%%%%%%%%%%%%%%%%%%%%%%%%%%%%%%%%%%%%%%%%%%%%%%%%%%%%%%%%%%%%%%%%%%%%%%%%%%%%%%%%%%%%%%%%%%%%%%
\section{Final comments}
%%%%%%%%%%%%%%%%%%%%%%%%%%%%%%%%%%%%%%%%%%%%%%%%%%%%%%%%%%%%%%%%%%%%%%%%%%%%%%%%%%%%%%%%%%%%%%%%
%

We have calculated the average of products of the characteristic moments 
$\langle \prod_{j=1}^m \det(I+x_j M) \rangle$.
We would also like to calculate the average of the quotient 
$$
\left\langle \frac{\prod_{j=1}^m \det(I+x_j M)}{\prod_{i=1}^l \det(I+y_i M)}
\right\rangle_n^{k_1,k_2,k_3}.
$$
Expressions for these quotients have been obtained for the classical groups (i.e., 
$(k_1,k_2,k_3)=(1,0,1), (0,1,1), (0,0,1)$ in our notation)   
in \cite{BG}, 
but the derivation of expressions for other cases remains an open problem.

\medskip
\medskip

\noindent
\textsc{Acknowledgements.}
The author would like to thank Professor Masato Wakayama for bringing to the author's attention  
the paper \cite{Mimachi}.

%<<<<<<<<<<<<<<<<<<<<<<<<<<<<<<<<<<<<<<<<<

\noindent
\textsc{Sho MATSUMOTO}\\
Faculty of Mathematics, Kyushu University.\\
Hakozaki Higashi-ku, Fukuoka, 812-8581 JAPAN.\\
\texttt{shom@math.kyushu-u.ac.jp}\\

\end{document}